\documentclass[leqno,12pt]{article}
\usepackage{latexsym}
 \usepackage{graphicx}
\usepackage{amsmath}
\usepackage{amssymb}
\usepackage{mathtools}
\usepackage{cite}
\usepackage{color}
\usepackage{subfig}
\usepackage[margin=1.3in, top=1.3in, bottom=1.3in]{geometry}
\usepackage{algorithm}
\usepackage{algorithmic}

\newcommand{\nc}{\newcommand}
\nc{\nt}{\newtheorem}
\nt{thm}{Theorem}[section]
\nt{cor}[thm]{Corollary}
\nt{prop}[thm]{Proposition}
\nt{lem}[thm]{Lemma}
\nt{defn}[thm]{Definition}
\nt{rem}[thm]{Remark}
\nt{exa}[thm]{Example}
\nt{ass}[thm]{Assumption}
\nt{alg}[thm]{Algorithm}
\nt{conj}[thm]{Conjecture}
\nt{claim}[thm]{Claim}
\nt{oracle}[thm]{Oracle}
\nc{\ip}[2]{\mbox{$\langle #1,#2 \rangle$}}
\nc{\pf}{\noindent{\bf Proof\ \ }}
\nc{\finpf}{\hfill{$\Box$}\linespace}
\nc{\linespace}{\vspace
{\baselineskip} \noindent}
\nc{\R}{{\mathbf R}}
\nc{\X}{{\mathbf X}}
\nc{\oR}{\overline{\R}}
\nc{\M}{\mathcal M}

\nc{\Rn}{{\mathbf R}^n}
\nc{\inT}{\mbox{\rm int}\,}
\nc{\cl}{\mbox{\rm cl}\,}

\def\tto{\;{\lower 1pt \hbox{$\rightarrow$}}\kern -12pt
           \hbox{\raise 2.8pt \hbox{$\rightarrow$}}\;}
\newenvironment{myequation}{\setcounter{equation}{\value{thm}}
   \begin{equation}}{\addtocounter{thm}{1}\end{equation}}

\nc{\bmye}{\begin{myequation}}
\nc{\emye}{\end{myequation}}

\begin{document}
\title{
The cost of nonconvexity in deterministic nonsmooth optimization
}
\author{
Siyu Kong
\thanks{ORIE, Cornell University, Ithaca, NY.
\texttt{sk3333@cornell.edu} 
}
\and
A.S. Lewis
\thanks{ORIE, Cornell University, Ithaca, NY.
\texttt{people.orie.cornell.edu/aslewis} 
\hspace{2cm} \mbox{~}
Research supported in part by National Science Foundation Grant DMS-2006990.}
}
\date{\today}
\maketitle

\begin{abstract}
We study the impact of nonconvexity on the complexity of nonsmooth optimization, emphasizing objectives such as piecewise linear functions, which may not be weakly convex.   We focus on a dimension-independent  analysis, slightly modifying a black-box algorithm of Zhang et al.\ \cite{mit} that approximates an $\epsilon$-stationary point of any directionally differentiable Lipschitz objective using $O(\epsilon^{-4})$ calls to a specialized subgradient oracle and a randomized line search.  Our simple black-box deterministic version, achieves $O(\epsilon^{-5})$ for any difference-of-convex objective, and $O(\epsilon^{-4})$ for the weakly convex case. Our complexity bound depends on a natural nonconvexity modulus, related, intriguingly, to the negative part of directional second derivatives of the objective, understood in the distributional sense.
\end{abstract}
\medskip

\noindent{\bf Key words:} nonsmooth optimization, nonconvex, Goldstein subgradient, \\
complexity, distributional derivative
\medskip

\noindent{\bf AMS Subject Classification:} 90C56, 49J52, 65Y20

\section{Introduction}
Minimizing nonsmooth nonconvex Lipschitz functions $f \colon \Rn \to \R$ is in general intractable.  A proxy approach, dating back to the 1970's \cite{gold_eps_stat}, can be viewed as seeking a point in $\Rn$ around which $f$ is differentiable on some cluster of nearby points at which the gradients have a small convex combination --- a small ``Goldstein subgradient'', in the language of \cite{gold_eps_stat}.  One algorithm accomplishing this goal, but without complexity analysis, appeared in \cite{mahdavi}.

A 2020 breakthrough \cite{mit} proposed an algorithm for this Goldstein subgradient problem with a complexity analysis depending on the radius $\delta$ of the cluster and the size $\epsilon$ of the subgradient, but {\em independent of the dimension} $n$.  The algorithm assumes directional differentiability of the objective $f$, and, relying on an associated directional subgradient oracle, uses an innovative randomized line search to achieve an efficiency guarantee of essentially $O(\epsilon^{-3}\delta^{-1})$.  Two subsequent papers \cite{somit,ddmit} point out that small random perturbations instead allow a standard subgradient oracle.

Two recent developments \cite{kornowski-shamir,jordan-lin-zampetakis}\footnote{We became aware of these concurrent independent works after completing the initial draft of this manuscript.} raise the question of {\em deterministic} algorithms for this problem.  While both prove positive results in the smooth case, and \cite{jordan-lin-zampetakis} thereby develops a ``white-box'' deterministic smoothing approach to the nonsmooth problem, both manuscripts also prove negative results for the general dimension-independent question.

Our contribution is to show, notwithstanding the negative results of \cite{kornowski-shamir,jordan-lin-zampetakis}, that a simple deterministic black-box version of the algorithm of \cite{mit} achieves, up to a {\em nonconvexity modulus} for the objective, a dimension-independent complexity of $O(\epsilon^{-4}\delta^{-1})$.  Our analysis covers interesting objectives such as piecewise linear functions, which are not even weakly convex.  We relate the nonconvexity modulus of the objective with its distributional second derivative, hinting at an intriguing relationship between such derivatives and algorithmic complexity in general optimization.

\section{The optimization problem and oracle}
Primarily to emphasize the elementary nature of our development, we adopt a rudimentary setting for our optimization problem.  On a real inner product space $\X$ with corresponding norm $|\cdot|$, we consider the problem of minimizing a function $f \colon \X \to \R$.  The objective $f$ may be neither smooth nor convex, and the space $X$ may be neither finite-dimensional, nor even complete.  

The method we develop, following \cite{mit}, relies on an underlying {\em directional subgradient map} $G \colon \X^2 \to \X$ associated with the objective $f$:  for all points $x \in \X$ and directions $e \in \X$, the G\^ateaux directional derivative
\[
f'(x;e) ~=~ \lim_{t \downarrow 0} \frac{1}{t}\big(f(x+te) - f(x)\big)
\]
exists and satisfies
\[
\ip{G(x,e)}{e} ~=~ f'(x;e).
\]
We say that $G$ is {\em $L$-bounded} for some constant $L>0$ if its norm $|G(x,e)|$ is never larger than $L$.
In applications, the objective function $f$ is $L$-Lipschitz and the vector $G(x,e)$ is a subgradient of some kind for $f$ at the point $x$ associated with the direction $e$, so we loosely refer to $G(x,e)$ as a ``subgradient''.  Nonetheless, we choose this rudimentary setting to emphasize again the elementary nature of our development, which makes no recourse to variational or Lipschitz analysis.

\begin{exa}[differentiable functions]
{\rm
If the function $f$ is $L$-Lipschitz and has a G\^ateaux derivative $\nabla f(x)$ at every point $x \in \X$, then the equation 
\[
G(x,e) ~=~ \nabla f(x)
\]
defines an $L$-bounded directional subgradient map.
}
\end{exa}

\begin{exa}[convex functions]
{\rm
For an $L$-Lipschitz convex function $f$ with convex subdifferential $\partial f$, any map $G$ satisfying 
\[
G(x,e) ~\in~ \mbox{argmax}\{\ip{g}{e} : g \in \partial f(x) \} \qquad \mbox{for all}~ x,e \in \X
\]
is an $L$-bounded directional subgradient map.
}
\end{exa}

More generally, we have the following example, which covers many interesting objectives, including the weakly convex case.  For a Lipschitz function $f \colon \Rn \to \R$, the {\em Clarke subdifferential} 
$\partial^c f(x)$ is the convex hull of the set of all limits of the form 
$\lim \nabla f(x^r)$ for points $x^r \to x$ in $\Rn$.  The function $f$ is {\em subdifferentially regular} when its G\^ateaux directional derivative satisfies
\[
f'(x;e) ~=~ \max\{\ip{g}{e} : g \in \partial^c f(x) \} \qquad \mbox{for all}~ x,e \in \Rn .
\]

\begin{exa}[subdifferentially regular functions]
{\rm
For an $L$-Lipschitz subdifferentially regular function $f \colon \Rn \to \R$, any map $G$ satisfying 
\[
G(x,e) ~\in~ \mbox{argmax}\{\ip{g}{e} : g \in \partial^c f(x) \} \qquad \mbox{for all}~ x,e \in \X
\]
is an $L$-bounded directional subgradient map.
}
\end{exa}

Notwithstanding the generality of this example, we emphasize that our framework is not restricted to objectives that are subdifferentially regular, as the following result shows.

\begin{prop}[directional Clarke subgradient maps] \label{clarke}
Any locally Lipschitz function $f \colon \Rn \to \R$ that is directionally differentiable has a directional subgradient map $G \colon \Rn \times \Rn \to \Rn$ satisfying $G(x,e) \in \partial^c f(x)$ for all $x,e \in \Rn$.
\end{prop} 

\pf
We just need to prove that if $f$ has a G\^ateaux directional derivative at the point $x \in \Rn$ in the direction $e \in \Rn$, then there exists a Clarke subgradient $g \in \partial^c f(x)$ satisfying $\ip{g}{e} = f'(x;e)$.  For $r=1,2,3,\ldots$, by the nonsmooth mean value theorem, there exists a point 
$x_r \in [x,x+\frac{1}{r}e]$ and a subgradient $g_r \in \partial^c f(x_r)$ satisfying
\[
f\Big(x+\frac{1}{r}e\Big) - f(x) ~=~ \ip{g_r}{\frac{1}{r} e}.
\]
Since the subdifferential $\partial^c f$ mapping is closed and locally bounded, any limit point $g$ of the sequence $\{g_r\}$ has the desired property.
\finpf

Minimizing nonconvex functions is in general intractable.  Instead we seek a point $x \in \X$ that is, in some sense, approximately critical.  To this end, we fix a constant $\delta > 0$, and consider what we call {\em Goldstein subgradients} at $x$:  vectors of the form
\[
\sum_{i=1}^k \lambda_i G(x_i,e_i)
\]
for a positive integer $k$, positive weights $\lambda_i$ summing to $1$, points 
$x_i$ in $B_{\delta}(x)$ (the closed ball of radius $\delta$ and center $x$) and directions $e_i \in \X$ for $i=1,2,\ldots,k$.  We denote the set of all Goldstein subgradients by $\partial_\delta f(x)$:  loosely, it consists of all convex combinations of subgradients at nearby points.  Strictly speaking, our notion is potentially smaller than the standard definition for Lipschitz $f \colon \Rn \to \R$, namely the closed convex hull of the set $\partial^c f\big(B_\delta(x)\big)$.

We can now state our goal, which relies on a second constant $\epsilon > 0$.
\bigskip

\noindent
{\bf Aim:}  Find a point $x \in \X$ and a Goldstein subgradient $g \in \partial_\delta f(x)$ such that $|g| \le \epsilon$.
\smallskip

\noindent
The development of \cite{mit} accomplishes this goal, explicitly in the setting of Proposition \ref{clarke}), and assuming that $f$ is directionally differentiable in the (stronger) Hadamard sense.  It relies on the following oracle.

\begin{oracle}[directional subgradient] \label{oracle}
{\rm
\mbox{} \\
{\bf Input:} 
\begin{itemize}
\item
a point $x \in \X$
\item
a direction $e \in \X$.
\end{itemize}
{\bf Output:}
\begin{itemize}
\item
the objective value $f(x)$
\item
the directional derivative $f'(x;e)$
\item
a subgradient-like vector $G(x,e)$.
\end{itemize}
}
\end{oracle}

In this work we rely on the same directional subgradient oracle, which we should emphasize is stronger than the standard subgradient oracle typically assumed.  In generic practice, for a Lipschitz objective $f$, we may expect that the algorithm never encounters points $x$ where $f$ is nondifferentiable, in which case any subgradient oracle simply returns $g = \nabla f(x)$.  More formally, however, we must consider nonsmooth points.  Undeterred, the authors of \cite{mit} argue that the stronger oracle is often a reasonable assumption.  However, our aim here is a fully deterministic algorithm. Accordingly, we use this same stronger oracle, but instead use a deterministic line search, much as in \cite{kornowski-shamir}.  To ensure termination, we make a mild assumption about the the directional behavior of the objective function $f$, similar in spirit to the idea of {\em semismoothness} \cite{mifflin-semismooth} common in nonsmooth computation, but simpler and weaker.  

We call the objective $f$ {\em directionally semismooth} if all points $x \in \X$ and directions $e \in \X$ satisfy
\[
\lim_{t \downarrow 0} f'(x+te;e) ~=~ f'(x;e).
\]
For comparison, following the definition in \cite{fin_2}, a Lipschitz function $f \colon \Rn \to \R$ is {\em semismooth} if it is directionally differentiable and for all points $x \in \Rn$ the directional derivative satisfies the stronger property
\[
f'(x+d;d) - f'(x;d) ~=~ o(d) \quad \mbox{as}~ d \to 0.
\]
Most Lipschitz functions in practice are semismooth, including in particular difference-of-convex functions and semi-algebraic functions \cite{semismooth}.  As we shall see, directional semismoothness suffices to guarantee termination of our algorithm, but before describing it, we focus first on the line search.

\section{A simple line search}
We pose the line search problem as a self-contained question.
Consider points $p<q$ in $\R$ and a function $h \colon [p,q] \to \R$ satisfying 
$h(p)>h(q)$.  Suppose that $h$ is right-differentiable on the interval $[p,q)$, and an oracle returns, for any input $t \in [p,q)$, the value $h(t)$ and the right derivative
\[
h'_+(t) ~=~ \lim_{s \downarrow t} \frac{h(s) - h(t)}{s-t}
\]
(possibly extended-valued).
How difficult is it to find a point $t$ satisfying $h'_+(t) < 0$? 

When $h$ is Lipschitz, the most basic randomized strategy --- uniformly sampling random points $t$ in the interval --- solves this problem with high probability.  Denoting the Lipschitz constant by $L$, the right derivative $h'_+$ always lies in the interval $[-L,L]$.  Denote the measure of the set where $h'_+ \ge 0$ by $\lambda$.  Then, providing the average slope satisfies
\[
\frac{h(q)-h(p)}{q-p} ~=~ -\sigma ~<~ 0,
\]
the fundamental theorem of calculus implies
\[
-\sigma ~=~ \frac{1}{q-p} \int_p^q h'_+ ~\ge~ 
L\Big(\frac{\lambda}{q-p} - 1\Big).
\]
Hence the probability that a uniformly distributed random point $t \in [p,q]$ fails to satisfy $h'(t) < 0$ is no larger than $1-\frac{\sigma}{L}$.  Thus for small $\sigma$, using at least $\frac{L}{\sigma}$ independent samples, the probability of success is at least $\frac{1}{2}$.   

However, we seek a deterministic algorithm, so we instead consider the following simple method, similar in spirit to one described in \cite{ddmit}.  We repeatedly bisect the interval $[p,q]$, each time discarding the subinterval over which the function $h$ decreases the least. The algorithm checks whether the right derivative at the midpoint of the current interval is negative, in which case it terminates.

\begin{alg}[bisection]
\label{bisection}
{\rm
\begin{algorithmic}
\STATE
\IF{$h'_+(p) < 0$}	
\RETURN{$p$}
\ENDIF
\STATE	$l=p$
\STATE	$r=q$
\WHILE{not done}
\STATE  $m = \frac{1}{2}(l+r)$
\IF{$h'_+(m) < 0$}	
\RETURN{$m$}
\ELSIF{$	2h(m) < f(l) + f(r)$}
\STATE	$r=m$
\ELSE
\STATE	$l=m$
\ENDIF
\ENDWHILE
\end{algorithmic}
}
\end{alg}

\noindent
Notice that the algorithm initially calls the oracle at the left endpoint $p$ of the given interval, calculates the function value at the right endpoint $q$, and then calls the oracle once during each bisection. 

In general, this algorithm may fail to terminate.  It is easy to construct a Lipshitz function $h$ satisfying $h(p) > h(q)$ and yet the derivative of $h$ at the initial endpoint $p$ and at every midpoint $m$ is positive.  To rule out such oscillatory examples, we can rely on directional semismoothness of $h$, which in this univariate setting simply means that the right derivative exists and is right continuous,
\[
\lim_{t \downarrow \bar t} h'_+(t) ~=~ h'_+(\bar t),
\]
and that the left derivative also exists and is left continuous:
\[
h'_-(t) ~=~ \lim_{s \uparrow t} \frac{h(s) - h(t)}{s-t}
\qquad 
\mbox{satisfies}
\qquad
\lim_{t \uparrow \bar t} h'_-(t) ~=~ h'_-(\bar t).
\]

For Lipschitz functions $h$, these two properties amount exactly to {\em semismoothness}.  Most Lipschitz functions in practice are semismooth, including in particular, convex and concave functions, and piecewise smooth functions.  Furthermore, any linear combination of semismooth functions is semismooth.  When the function $h$ is semismooth, the Clarke subdifferential is given by 
\[
\partial^c h(t) ~=~ \mbox{conv}\{h'_-(t)\, ,\, h'_+(t)\},
\]
and the following property also holds \cite[Lemma 2]{mifflin-semismooth} and \cite[Lemma 2.2]{henrion-outrata}:
\[
\lim_{t \downarrow \bar t} h'_-(t) ~=~ h'_+(\bar t)
\qquad 
\mbox{and}
\qquad
\lim_{t \uparrow \bar t} h'_+(t) ~=~ h'_-(\bar t).
\]

Semismoothness is more than enough to prove termination of the line search.  The simple argument below also applies to nonlipschitz functions.

\begin{prop}[Termination of the line search]~~
\label{terminates}
Suppose that the function $h \colon [p,q] \to \R$ satisfies 
$h(p) > h(q)$, and that its left and right derivatives satisfy the semismoothness conditions
\begin{eqnarray*}
\lim_{t \downarrow \bar t} h'_+(t) ~=~ h'_+(\bar t) \qquad \mbox{for}~ \bar t \in [p,q) \\
\lim_{t \uparrow \bar t} h'_+(t) ~=~ h'_-(\bar t) \qquad \mbox{for}~ \bar t \in (p,q].
\end{eqnarray*}
Then Algorithm~\ref{bisection} terminates.
\end{prop}

\pf
If the algorithm does not terminate, then it generates monotonic sequences $l_k \uparrow$ and $r_k \downarrow$, satisfying $r_k-l_k \to 0$,
\[
h'_+(l_k) ~\ge~ 0, \qquad \mbox{and} \qquad 
\frac{h(r_k) - h(l_k)}{r_k-l_k} ~\le~ -\sigma ~<~ 0
\]
for each iteration $k=0,1,2,\ldots$.  Denote the two sequences' mutual limit by $\bar m$.  Semismoothness ensures $h'_+(l_k) \to h'_-(\bar m)$, so $h'_-(\bar m) \ge 0$.

If $r_k = \bar m$ for all large $k$, then
\[
\frac{h(r_k) - h(l_k)}{r_k-l_k} ~=~ \frac{h(\bar m) - h(l_k)}{\bar m -l_k} ~\to~ h'_-(\bar m),
\]
which is a contradiction.  Hence for all large $k$ we have $q > r_k > \bar m$, so $h'_+(r_k) \ge 0$, and
semismoothness ensures
$h'_+(r_k) \to~ h'_+(\bar m)$, so $h'_+(\bar m) \ge 0$.  We deduce
\[
h(\bar m) - h(l_k) ~\ge~ -\frac{\sigma}{2}(\bar m - l_k)
\qquad \mbox{and} \qquad h(r_k) - h(\bar m) ~>~ -\frac{\sigma}{2}(r_k - \bar m).
\]
Adding now gives a contradiction.
\finpf

\section{The optimization algorithm}
To minimize a locally Lipschitz function $f \colon \X \to \R$ using the directional subgradient oracle described above, we study the following algorithm.  The method we describe is essentially that of \cite{mit}, but with the deterministic line search described in the preceding section.  Given a current point $x \in \X$ and a current Goldstein subgradient,
the method proceeds as follows.

\begin{alg}[nonsmooth minimization]
\label{minimization}
{\rm
\begin{algorithmic}
\STATE
\WHILE{not done}	
		\IF{$|g| \le \epsilon$}
			\RETURN{$x$}								\hfill \COMMENT{{\em Small subgradient so stop.}}
		\ENDIF
		\STATE $\hat g = \frac{g}{|g|}$					\hfill \COMMENT{{\em Normalize subgradient.}}
		\STATE  $x' = x - \delta \hat g$				\hfill \COMMENT{{\em Trial step of fixed length}.}
		\IF{$f(x) - f(x') \ge \frac{\delta\epsilon}{3}$}	
			\STATE $x=x'$								\hfill \COMMENT{{\em Sufficient decrease so update point.}}
			\STATE $g = G(x,0)$			\hfill \COMMENT{{\em Re-initialize subgradient.}}
		\ELSE
			\STATE										\hfill \COMMENT{{\em Insufficient decrease so update subgradient.}}
			\STATE Define $h$ on $[0,\delta]$ by
			\begin{eqnarray*}
			x(t) &=& x + (t-\delta) \hat g \\
			h(t) & = & f\big( x(t) \big) ~-~ \frac{\epsilon t}{2}. 
			\end{eqnarray*}
			Apply Algorithm \ref{bisection} (bisection), using the formula
			\[
			h'_+(t) ~=~ \ip{G(x(t),\hat g)}{\hat g} ~-~ \frac{\epsilon}{2}
			\]
			to find $t \in [0,\delta]$ satisfying
			$h'_+(t) < 0$.
			\STATE $g = \mbox{shortest vector in}~ [g,G(x(t),\hat g)]$	
		\ENDIF
\ENDWHILE
\end{algorithmic}
}
\end{alg} 

When the objective $f$ is directionally semismooth, Algorithm \ref{bisection} terminates by Proposition \ref{terminates}, which in turn guarantees termination of Algorithm \ref{minimization}, as we shall now prove.  We use the following simple tool, following \cite{mit}.

\begin{lem}
If two vectors $g,g' \in B_L(0)$ satisfy $\ip{g'}{g} \le \frac{1}{2}|g|^2$, then the shortest vector $g''$ in the line segment $[g,g']$ satisfies 
\[
|g''|^2 ~\le~ |g|^2 \Big(1 - \frac{|g|^2}{16L^2} \Big).
\]
\end{lem}

\pf
For all $t \in [0,1]$ we have
\[
|g''|^2 
~\le~ 
|g + t(g'-g)|^2 
~=~ 
|g|^2 + t^2|g'-g|^2 + 2t \ip{g}{g'-g}
~ \le ~
(1-t)|g|^2 + 4L^2t^2.
\]
Setting $t = \frac{|g|^2}{8L^2}$ proves the result.
\finpf

We can now prove the validity of the algorithm, again imitating parts of the argument in \cite{mit}, which we reproduce for ease of reading.

\begin{thm}[finite termination] \label{finite}
Suppose that we apply Algorithm \ref{minimization} to a directionally semismooth function 
$f \colon \X \to \R$ that is bounded below.  Suppose that the directional subgradient map in Oracle \ref{oracle} is $L$-bounded.  Then the algorithm terminates with a point $x \in \X$ and a subgradient $g \in \partial_\delta f(x)$ satisfying\ $|g| \le \epsilon$.
\end{thm}

\pf
Suppose that the current subgradient $g$ is {\em inadequate}, in the sense that, in the terminology of the algorithm description, it is neither small, nor generates sufficient decrease.
We then apply the bisection method, Algorithm \ref{bisection}, to the given function $h$.  The initial interval is $[p,q] = [0,\delta]$, and the average slope satisfies
\[
\frac{h(q)-h(p)}{q-p} ~=~ 
\frac{1}{\delta} \Big( f(x) - \frac{\delta\epsilon}{2} - f(x')\Big) ~<~ - \frac{\epsilon}{6}.
\]
We thus arrive at a subgradient $g' \in \partial_{\delta} f(x)$ satisfying 
$\ip{g'}{\hat g} < \frac{\epsilon}{2}$.  We deduce $\ip{g'}{g} < \frac{|g|^2}{2}$.  The algorithm replaces the current subgradient $g$ by the shortest vector $g''$ in the line segment $[g,g']$.  
Since $g'' \in \partial_\delta f(x)$, we can repeat this shortening process, providing that $g''$ is also inadequate.  Suppose that the subgradient $g$ remains inadequate after completing $k$ such shortening steps.  Let 
$\rho_i$ denote the quantity $\frac{|g|^2}{16L^2}$ after $i=0,1,2,\ldots,k$ steps.  
Then $\rho_0 \le \frac{1}{16}$ and for each $i$ we have
\[
0 ~<~ \rho_{i+1} ~\le~ \rho_i(1-\rho_i),
\]
so
\[
\frac{1}{\rho_{i+1}} ~\ge~ \frac{1}{\rho_i} + \frac{1}{1 - \rho_i} ~>~ \frac{1}{\rho_i} + 1.
\]
Consequently,
\[
\frac{1}{\rho_k} ~\ge~ 16 + k,
\]
so we deduce
\[
\frac{\epsilon^2}{16L^2} ~<~ \frac{|g|^2}{16L^2} ~\le~ \frac{1}{16+k}.
\]
Hence after no more than
\[
16 \Big(\frac{L^2}{\epsilon^2} - 1 \Big)
\]
shortening steps, each requiring one line search, we arrive at an adequate subgradient 
$g \in \partial_\delta f(x)$.  To summarize, starting at any point with an inadequate subgradient, we require no more than $\frac{16L^2}{\epsilon^2}$ line searches before finding an adequate subgradient $g$.

There are now two possibilities.  Either the subgradient $g$ satisfies $|g| \le \epsilon$, in which case we stop, or we perform a reduction step, replacing the current point $x$ by $x - \delta \frac{g}{|g|}$, thereby decreasing the objective value by at least the quantity $\frac{\delta\epsilon}{3}$.  Since the objective is bounded below, beginning from the initial point $x_0$, this procedure terminates after no more than 
$\lceil \frac{3}{\delta\epsilon}\big(f(x_0) - \min f)\big)\rceil$ reduction steps
and hence
\bmye \label{searches}
\left\lceil \frac{3(f(x_0) - \min f)}{\delta\epsilon} \right\rceil \cdot \frac{16L^2}{\epsilon^2} 
\emye
line searches.  Proposition \ref{terminates} ensures that each line search requires only finitely many oracle calls, completing the proof.
\finpf

\section{Complexity of the line search}
To complete our complexity analysis for the minimization algorithm, we simply need to bound the number of oracle calls needed by each line search, and multiply by our bound (\ref{searches}) on the number of line searches.  Consider, therefore, the bisection method.  When the function $h \colon [p,q] \to \R$ is convex, the problem is trivial: since
\[
h(p) > h(q) \quad \Rightarrow \quad h'_+(p) < 0,
\]
the algorithm terminates at the first oracle call.  More generally, we proceed by correcting any lack of convexity in $h$ by adding a convex perturbation $s \colon [p,q] \to \R$.  

Recall that, for any interval $J$, a function $h \colon J \to \R$ is {\em difference-of-convex} when there exists a convex function $s \colon J \to \R$ such that $h+s$ is also convex.  For such functions, we have the following tool.

\begin{lem}
\label{estimate}
Consider a function $h \colon [p,q] \to \R$, and a convex function $s \colon [p,q] \to \R$ such that $h+s$ is also convex.
For any points $x < y$ in the interval $[p,q]$, if $h'_+(x) \ge 0$, then
\[
\frac{h(y) - h(x)}{y-x} ~\ge~ s'_+(x)-s'_-(y).
\]
\end{lem}

\pf
Denote the convex function $h+s$ by $r$.
The convex functions $r$ and $s$ satisfy
\[
r'_+(x) ~\in~ \partial r(x) \qquad \mbox{and} \qquad
s'_-(y) ~\in~ \partial s(y).
\]
Hence
\begin{eqnarray*}
s'_-(y)(x-y) 
&\le& 
s(x)-s(y) ~=~ r(x)-h(x)-r(y)+h(y) \\ 
&\le&
h(y)-h(x) + r'_+(x)(x-y) ~\le~  h(y)-h(x) + s'_+(x)(x-y),
\end{eqnarray*}
and the result follows.  
\finpf

We can then use the change in derivative of the necessary perturbation $s$ to bound the number of iterations in the line search.

\begin{thm} \label{time}
Consider a function $h \colon [p,q] \to \R$, and a convex function $s \colon [p,q] \to \R$ such that $h+s$ is also convex.  If the bisection method, Algorithm \ref{bisection} 
evaluates $h'_+$, the right-derivative, $k \ge 1$ times without terminating, then
\[
\frac{h(q) - h(p)}{q-p} ~\ge~ \frac{s'_+(p)-s'_-(q)}{k}.
\]
\end{thm}

\pf
We proceed by induction on the number of evaluations $k=1,2,3,\ldots$.  The case $k=1$ follows immediately from Lemma \ref{estimate}, by setting $x=p$ and $y=q$.  

Suppose that the result holds for any points $p<q$ and for any number of evaluations no larger than $k$.  Now consider an instance of the algorithm that completes $k+1$ evaluations.  After the first evaluation, consider the midpoint $m = \frac{1}{2}(p+q)$.  There are two possible cases, depending on whether or not 
\bmye \label{dichotomy}
2h(m) ~<~ h(p)+h(q).
\emye
We consider them in turn.

Suppose first that inequality (\ref{dichotomy}) holds.  After the first evaluation, the algorithm makes $k$ further evaluations, beginning with the initial interval $[p,m]$.  Hence, by the induction hypothesis,
\[
\frac{h(m) - h(p)}{m-p} ~\ge~ \frac{s'_+(p)-s'_-(m)}{k}.
\]
Since the algorithm did not terminate during the first two evaluations, we know  
$h'_+(m) \ge 0$.  Applying Lemma \ref{estimate} with $x=m$ and $y=q$ shows
\[
\frac{h(q) - h(m)}{q-m} ~\ge~ s'_+(m)-s'_-(q).
\]
Hence
\begin{eqnarray*}
s'_+(p)-s'_-(q) 
&\le& \big(s'_+(p)-s'_-(m)\big) ~+~ \big(s'_+(m)-s'_-(q)\big) \\
&\le&
k \frac{h(m) - h(p)}{m-p} ~+~ \frac{h(q) - h(m)}{q-m} \\
&=&
(k-1) \frac{h(m) - h(p)}{m-p} ~+~ \Big( \frac{h(m) - h(p)}{m-p} ~+~ \frac{h(q) - h(m)}{q-m} \Big) \\
&=&
(k-1) \frac{h(m) - h(p)}{m-p} ~+~  2\frac{h(q) - h(p)}{q-p} \\
&<&
(k-1) \frac{h(q) - h(p)}{q-p} ~+~  2\frac{h(q) - h(p)}{q-p} ~=~ (k+1) \frac{h(q) - h(p)}{q-p},
\end{eqnarray*}
as required.

The case where inequality (\ref{dichotomy}) fails is similar.  After the first bisection, the algorithm makes $k$ further bisections, beginning with the initial interval $[m,q]$.  Hence, by the induction hypothesis,
\[
\frac{h(q) - h(m)}{q-m} ~\ge~ \frac{s'_+(m)-s'_-(q)}{k}.
\]
Since the algorithm did not terminate during the first bisection, we know  
$h'_+(p) \ge 0$.  Applying Lemma \ref{estimate} with $x=p$ and $y=m$ shows
\[
\frac{h(m) - h(p)}{m-p} ~\ge~ s'_+(p)-s'_-(m).
\]
Hence
\begin{eqnarray*}
s'_+(p)-s'_-(q) 
&\le& \big(s'_+(p)-s'_-(m)\big) ~+~ \big(s'_+(m)-s'_-(q)\big) \\
&\le&
\frac{h(m) - h(p)}{m-p} ~+~ k\frac{h(q) - h(m)}{q-m} \\
&=&
\Big( \frac{h(m) - h(p)}{m-p} ~+~ \frac{h(q) - h(m)}{q-m} \Big)
 ~+~ (k-1) \frac{h(q) - h(m)}{q-m} \\
&=&
2\frac{h(q) - h(p)}{q-p} ~+~ (k-1) \frac{h(q) - h(m)}{q-m}  \\
&\le&
2\frac{h(q) - h(p)}{q-p} ~+~ (k-1) \frac{h(q) - h(p)}{q-p} ~=~ (k+1) \frac{h(q) - h(p)}{q-p},
\end{eqnarray*}
as required.
\finpf

\begin{defn}
{\rm
Given any interval $J \subset \R$, the {\em concave deviation} of a function $h \colon J \to \R$ is the infimum of the Lipschitz constants of convex functions $s \colon J \to \R$ such that the sum $h+s$ is also convex.  
}
\end{defn}
Consider, for example, a $\rho$-{\em weakly convex} function $h$, for some constant 
$\rho  \ge 0$, meaning that the function 
$t \mapsto h(t) + \frac{\rho}{2}t^2$ is convex.  

\begin{prop}
\label{weakly}
Any \mbox{$\rho$-weakly} convex function $h \colon [p,q] \to \R$, for $\rho \ge 0$, has concave deviation at most $\frac{\rho}{2}(q-p)$.
\end{prop}

\pf
The function $s(t) = \frac{\rho}{2}(t-\frac{p+q}{2})^2$ is convex, with Lipschitz constant $\frac{\rho}{2}(q-p)$, and $h+s$ is also convex.
\finpf

The concave deviation for functions that are not weakly convex may shrink more slowly than the length of the interval.  For example, on the interval $[-\delta,\delta]$, the function $h(t) = -|t|^\frac{3}{2}$  has concave deviation $\frac{3}{2}\sqrt{\delta}$, and the piecewise linear function $-|\cdot|$ has concave deviation $1$. 

The concave deviation of a function $h \colon [p,q] \to \R$ that is difference-of-convex may not be finite.  An example is the concave function $\sqrt{\cdot}$ on the interval $[0,1]$.  However, if $h$ extends to a difference-of-convex function on an open interval containing the interval $[p,q]$, then its nonconvexity bound on $[p,q]$ must be finite, since we can write $h$ as a difference of convex functions, each of which must be Lipschitz on $[p,q]$.

Consider any continuous semi-algebraic function $h \colon \R \to \R$.  We can partition $\R$ into finitely-many closed intervals $J$ such that each restriction $h|_J$ is either convex or concave.  In general, such a partition may not guarantee that $h$ is difference-of-convex: an example is the function $x^{\frac{1}{3}}$.  However, if $h$ is also Lipschitz, then each convex or concave ingredient $h|_J$ extends to a corresponding convex or concave Lipschitz function on $\R$, and from these we can easily decompose $h$ into a difference of convex Lipschitz functions.  Thus all semi-algebraic functions on $\R$ are difference-of-convex, with finite concave deviation on any bounded interval.

\begin{cor}[Line search complexity]~
\label{line}
If a function $h \colon [p,q] \to \R$ has finite concave deviation $M$ and average rate of decrease
\[
\sigma ~=~ -\frac{h(q)-h(p)}{q-p}  ~>~ 0,
\]
then the number of evaluations of $h'_+$, the right-derivative, required before the bisection method, Algorithm \ref{bisection}, terminates is no more than $1 + \big\lfloor \frac{2M}{\sigma} \big\rfloor$.
\end{cor} 

\pf
Suppose that the bisection method evaluates the left-derivative
$k \ge 1$ times without terminating.  Fix any value $M' > M$.
By assumption, there exists a convex function $s$ with Lipschitz constant less than $M'$ such that the sum $h+s$ is also convex.  From Theorem \ref{time}, we deduce the inequalities
\[
-\sigma ~\ge~ \frac{s'_+(p)-s'_-(q)}{k} ~>~ \frac{-2M'}{k}
\]
so $k < \frac{2M'}{\sigma}$.  Since $M'$ was arbitrary, we deduce $k \le \frac{2M}{\sigma}$,
and hence $k \le \lfloor \frac{2M}{\sigma} \rfloor$.  The result follows.
\finpf

Given an open interval $I$, consider a difference-of-convex function $h \colon I \to \R$.  As observed in \cite{hartman}, such functions are characterized by having left and right derivatives everywhere, which furthermore are of bounded variation on every compact interval in $I$.  Any such function $h$ also has a second derivative $D^2 h$ in the distributional sense:  in general it is a signed Radon measure on the interval $I$ (see \cite{dudley}).  

We review briefly the underlying construction.  Consider any convex function  $s \colon I \to \R$.  Its right derivative $s'_+$ is nondecreasing and right-continuous, and hence defines a nonnegative Radon measure $D^2 s$ on the interval $I$ via the property
\[
(D^2 s)(p,q] ~=~ s'_+(q) - s'_+(p) \qquad \mbox{for all}~ p<q ~\mbox{in}~ I.
\]
More generally, for any difference-of-convex function $h \colon I \to \R$, consider any convex function  
$s \colon I \to \R$ such that $h+s$ is also convex.  The second derivative $D^2h$ is just the signed measure $D^2(h+s) - D^2 s$,
which is independent of the choice of $s$.  Convexity of $h$ is characterized by the property 
$D^2 h \ge 0$.  More generally, the Hahn decomposition decomposes $D^2 h$ uniquely into a difference of nonnegative Radon measures supported on disjoint sets.  These positive and negative parts satisfy  
\[
D^2 h ~=~ (D^2 h)^+ ~-~ (D^2 h)^-.
\]

\begin {thm}
If the function $h \colon [p,q] \to \R$ is difference-of-convex, then its concave deviation is 
\[
\frac{1}{2} (D^2 h)^- (p,q).
\]
\end{thm}

\pf
For any convex function $s : [p,q] \to \R$ such that $h+s$ is also convex, the second derivatives satisfy
\[
0 ~\le~ D^2 s ~=~ D^2 (h+s) - D^2 h ~\ge~ -D^2 h,
\]
so $D^2 s \ge (D^2 h)^-$.
If $s$ is $M$-Lipschitz on $[p,q]$, then 
\[
M 
~\ge~
\max \{ -s'_+(p) , s'_-(q) \} ~\ge~
\frac{s'_-(q) - s'_+(p)}{2} \\
~=~
\frac{1}{2} (D^2 s) (p,q) ~\ge~ \frac{1}{2} (D^2 h)^- (p,q).
\]

If the right-hand side is infinite, this completes the proof, so suppose that it is finite.  Define a function $g \colon (p,q] \to \R$ by
\[
g(t) ~=~ (D^2 h)^- (p,t).
\]
Then $g$ is a nonnegative nondecreasing left-continuous function that is bounded above and 
$g(t) \downarrow 0$ as $t \downarrow p$.  Now define a convex function $s \colon [p,q] \to \R$ by
\[
s(t) ~=~ \int_p^t g(\tau)\,d\tau.
\]
Then $s'_-(t) = g(t)$ for all $t \in (p,q]$ and $s'_+(p) = 0$.  Furthermore we have
\[
D^2(h+s) ~=~ D^2 h + D^2 s ~=~ D^2 h + (D^2 h)^- ~=~ (D^2 h)^+ ~\ge~ 0,
\]
so the sum $h+s$ is also convex.  The function
\[
t ~\mapsto~ \tilde s(t) ~=~ s(t) - \frac{1}{2}s'_-(q) t
\]
is also convex, as is $h+\tilde s$, and the function $\tilde s$ has Lipschitz constant
\[
-\tilde s'_+(p) ~=~ \tilde s'_-(q) ~=~ \frac{1}{2} s'_-(q) ~=~ \frac{1}{2} (D^2 h)^- (p,q).
\]
This completes the proof.
\finpf

As an illustration, we have the following result.

\begin{cor}[Piecewise linear functions] \label{pw}
Consider any continuous piecewise linear function $h \colon [p,q] \to \R$, with $m$ derivative discontinuities $t_1 < t_2 < \cdots < t_m$ in the interval $(p,q)$. Define $t_0 = p$ and $t_{m+1}=q$, and let $g_i$ be the value of the derivative on the  interval $(t_i,t_{i+1})$ for $0 \le i \le m$.  Then $h$ has concave deviation  
\[
\frac{1}{2}\sum_{i=1}^m (g_{i-1} - g_i)^+.
\]
If $h$ is $L$-Lipschitz, then this bound is no larger than $\big\lceil \frac{m}{2} \big\rceil L$.
\end{cor}

\pf
Denoting a unit point mass at the point $t$ by $\delta_t$, we have
\[
D^2 h ~=~ \sum_i (g_i - g_{i-1})\delta_{t_i},
\]
and hence 
\[
(D^2 h)^- ~=~ \sum_i (g_{i-1} - g_i)^+ \delta_{t_i},
\]
from which the claimed equation follows.  The inequality is an easy consequence, using $|g_i| \le L$ for each $i$.
\finpf

As an illustration, we present an example that underlines why the line search complexity estimate in Corollary \ref{line} is the best that we can expect in general.  In outline, while the functions $h \colon [p,q] \to \R$ that we consider satisfy $h(p) > h(q)$, their derivatives may often be positive.

\begin{exa}[Optimality of the line search]
{\rm
Consider any constant $M > 0$, and a set $T \subset [0,1)$ of cardinality strictly less than $2M$.  Then there exists a function $h \colon [0,1] \to \R$ with concave deviation less than $M$, that satisfies $h(0)=0$ and $h(1)=-1$, and has strictly positive derivative throughout $T$.

To see this, suppose first $T \subset (0,1)$.  (The case when $T$ contains zero is an easy modification.)  Enumerate the points in increasing order:
\[
t_1 < t_2 < \cdots < t_k,
\]
where $k<2M$.  Define $h(0) = 0$ and $h(1)=-1$.
Fix any small $\gamma > 0$, and define
\[
h(t_i - \gamma) ~=~ -t_i - \gamma^2 \quad \mbox{and} \quad 
h(t_i + \gamma) ~=~ -t_i + \gamma^2 \quad \mbox{for}~ i=1,2,\ldots,k.
\]
At intermediate points in $[0,1]$, define $h$ by linear interpolation.  A quick calculation, using Corollary \ref{pw}, shows that $h$ has concave deviation
\[ 
\frac{1}{2}
\sum_{i=1}^{k-1}
\Big(
\frac{t_{i+1} - t_i + 2\gamma^2}{t_{i+1} - t_i - 2\gamma} + \gamma
\Big)
~+~ 
\frac{1}{2}
\Big(
\frac{1 - t_k + \gamma^2}{1 - t_k - \gamma} + \gamma
\Big)
~=~ \frac{k}{2} + O(\gamma) ~<~ M
\]
providing that $\gamma$ is sufficiently small.

Now consider any line search method applicable to functions $h \colon [0,1] \to \R$ satisfying $h(0)=0$ and $h(1) = -1$, relying on evaluations of the value $h$ and the right derivative $h'_+$ at points chosen one-by-one, and terminating once a derivative is negative.  Suppose that the method is guaranteed to terminate after at most $k$ queries providing that the underlying function $h$ has concave deviation strictly less than some given value $M > 0$.  The example above proves $k \ge 2M$.
}
\end{exa} 

\section{Multivariate functions}
To understand the complexity of Algorithm \ref{minimization} (nonsmooth minimization), we apply our analysis in the previous section to restrictions of multivariate objectives $f$ to line segments. 
For any convex set $C \subset \X$, consider a function $f \colon C \to \R$.  Given any length 
$\delta > 0$, let $\Lambda(\delta)$ denote the supremum over all points $x,y \in C$ with 
$|x-y| \le \delta$ of the concave deviation for the function $h \colon [0,\delta] \to \R$ defined by
\bmye \label{1D}
h(t) ~=~ f\Big(x + \frac{t}{\delta}(y-x)\Big).
\emye
We call the function $\Lambda \colon \R_{++} \to [0,+\infty]$ the {\em nonconvexity modulus} for $f$.   The following illustration follows immediately from Proposition \ref{weakly}.

\begin{prop}
The nonconvexity modulus of any $\rho$-weakly convex function (for $\rho \ge 0$) satisfies
\[
\Lambda(\delta) ~\le~ \frac{\rho\delta}{2}.
\]
\end{prop}

More generally, the function $f$ is {\em difference-of-convex} when there exists a convex function $q \colon C \to \R$ such that $f+q$ is also convex.  

\begin{prop}
Consider a convex set $C \subset \X$ and functions $f,q \colon C \to \R$ with both $q$ and $f+q$ convex.  If $q$ is $M$-Lipschitz, then the nonconvexity modulus of $f$ satisfies $\Lambda(\delta) \le M$ for all 
$\delta > 0$.
\end{prop}

\pf
Consider any points $x,y \in C$ with $|x-y| \le \delta$, and the function $h$ defined by equation (\ref{1D}).  The function 
$s \colon [0,\delta] \to \R$ defined by
\[
s(t) ~=~ q\Big(x + \frac{t}{\delta}(y-x)\Big)
\]
is convex, and $M$-lipschitz, and $h+s$ is convex, so the concave deviation of $h$ is no larger than $M$.  The result follows.
\finpf

\noindent
As a consequence, we deduce the following result.

\begin{cor}
Consider any convex sets $C \subset C' \subset \Rn$, where $C$ is nonempty and compact and $C'$ is open, and any difference-of-convex function $f \colon C' \to \R$.  Then the nonconvexity modulus of the restriction $f_C$ is uniformly bounded: there exists a finite constant $M$ such that 
$\Lambda(\delta) \le M$ for all $\delta > 0$.
\end{cor}

\begin{cor}
If a function $f \colon \X \to \R$ is the difference $p-q$ between a convex function $p \colon \X \to \R$ and a polyhedral convex function $q \colon \X \to \R$, then the nonconvexity modulus of $f$ is no larger than any Lipschitz constant for $q$.
\end{cor}

More generally, consider a continuous function $f \colon \X \to \R$ that is {\em semi-linear}, in the sense that $\X$ is a finite union of polyhedra, on each of which the function $f$ is affine.  Any such function has a Lipschitz constant $L$, and furthermore a uniform upper bound $m$ on the number of possible gradient discontinuities in any function of the form (\ref{1D}).  By Corollary \ref{pw}, we deduce that the nonconvexity modulus 
$\Lambda(\delta)$ is no larger than $\big\lceil \frac{m}{2} \big\rceil L$.

Returning to our analysis of Algorithm \ref{minimization}, we are ready for our main result.

\begin{thm}[complexity of minimization]
Given two constants $\delta,\epsilon > 0$, consider a convex set $C \subset \X$, a function 
$f \colon C \to \R$ that is bounded below, an associated $L$-bounded directional subgradient map 
$G \colon \X^2 \to \X$, and an initial point $x_0 \in C$ such that 
\[
f(x) \le f(x_0) \quad \mbox{and} \quad |y| \le \delta \quad \Rightarrow \quad x+y \in C.
\]
Suppose that $f$ has finite nonconvexity modulus $\Lambda(\delta)$.  Then Algorithm \ref{minimization} (nonsmooth minimization) requires at most
\[
\left\lceil \frac{3(f(x_0) - \min f)}{\delta\epsilon}\right\rceil \cdot \frac{16L^2}{\epsilon^2} \cdot 
\Big(1 + \Big\lfloor \frac{12\Lambda(\delta)}{\epsilon} \Big\rfloor \Big)
\]
calls to Oracle \ref{oracle} (directional subgradient)  to find a point $x \in \X$ and a Goldstein subgradient $g \in \partial_\delta f(x)$ satisfying\ $|g| \le \epsilon$.
\end{thm}

\pf
Corollary~\ref{line} with
$\sigma = \frac{\epsilon}{6}$ shows that the bisection method requires at most
\[
1 + \Big\lfloor \frac{12\Lambda(\delta)}{\epsilon} \Big\rfloor
\]
oracle calls to terminate.  Then, one further call returns the desired subgradient
$g' \in \partial_\delta f(x)$ satisfying $\ip{g'}{g} < \frac{|g|^2}{2}$.  Multiplying by the bound (\ref{searches}) on the number of line searches completes the proof.
\finpf

\section{Appendix:  distributional second derivatives}
We saw previously that the concave deviation of a univariate function $h$ is determined by the negative part of its second distributional derivative $D^2h$.  In our application, we consider functions $h$ that are restrictions of the underlying objective $f$ to line segments of fixed length $\delta$.  We would therefore expect the nonconvexity modulus of $f$ to be related to its own distributional second derivative.  Here, we explore that relationship informally.

Consider a locally Lipschitz function $f \colon \Rn \to \R$.  The distributional derivative of $f$ is an $n$-vector $Df$, entries of which are distributions --- linear functionals on the space of smooth, compactly supported functions $g \colon \Rn \to \R$ ({\em test functions}), that are continuous with respect to uniform convergence on compact sets.  We can define $Df$ through the relationship
\[
\ip{u^T (Df)}{g} ~=~ - \int  f (u^T \nabla g)
\]
for all vectors $u \in \Rn$ and test functions $g \colon \Rn \to \R$.  However, by a suitable version of Rademacher's Theorem \cite[Section 6.2, Theorem 1]{evans-gariepy}, the gradient $\nabla f$ exists almost everywhere and is essentially bounded, and satisfies
\[
\ip{u^T (Df)}{g} ~=~ \int  g (u^T \nabla f).
\]
In standard terminology \cite{evans-gariepy}, we can identify the classical gradient $\nabla f$ with both the distributional derivative $Df$ and the ``weak'' derivative of $f$.  

The second distributional derivative of $f$ is an $n$-by-$n$ matrix $D^2 f$, entries of which are distributions.  We can define $D^2 f$ through the relationship
\[
\ip{u^T (D^2 f) v}{g} ~=~ - \int (u^T \nabla f) (v^T \nabla g)
\]
for all vectors $u,v \in \Rn$ and test functions $g$.  
If $f$ is smooth, then $D^2 f$ is just the matrix-valued measure with density $\nabla^2 f$.  More generally we must consider $D^2 f$ as a distribution, but at least for {\em convex} functions we can be more specific:  it is a positive-semidefinite-valued Radon measure (\cite{dudley} and 
\cite[Section 6.3]{evans-gariepy}).

\begin{exa}[A piecewise linear function]~~ \label{hinge}
{\rm
Consider the convex function $f \colon \R^2 \to \R$ defined by $f(x) = x_1^+$.  For any vectors $u,v \in \R^2$ and smooth, compactly supported function $g \colon \R^2 \to \R$, we have
\begin{eqnarray*}
\ip{u^T (D^2 f) v}{g} 
&=&
 - \int_{x_1 > 0} u_1 v^T \nabla g(x)\, dx
~=~ - u_1 \int_{x_1 > 0} \mbox{div}(g(x) v)\, dx \\
&=& -u_1 \int_{\R} (-e_1)^T (g(\bigl[ \begin{smallmatrix} 0 \\ y \end{smallmatrix} \bigr]) v)\, dy ~=~ u_1 v_1 \int g(\bigl[ \begin{smallmatrix} 0 \\ y \end{smallmatrix} \bigr])\, dy,
\end{eqnarray*}
by the Gauss-Green formula.
Thus $D^2 f$ is the matrix $\bigl[ \begin{smallmatrix} \mu & 0 \\ 0 & 0 \end{smallmatrix} \bigr]$, where the measure $\mu$ is related to Lebesgue measure $\lambda$ via
\bmye \label{measure}
\mu(S) ~=~ \lambda\{ s \in \R : \bigl[ \begin{smallmatrix} 0 \\ s \end{smallmatrix} \bigr] \in S \},
\emye
for all measurable subsets of $S \subset \R^2$.
}
\end{exa}

For a more general understanding, we begin with the univariate case.
\begin{exa}[Univariate convex functions] \label{univariate}~~~
{\rm
Consider a convex function $f \colon \R \to \R$.    
For $0 < \gamma < 1$, we can construct a smooth approximation $h_\gamma \colon \R \to [0,1]$ of the standard step function, with the following properties: 
\[
h_\gamma (t) ~=~
\left\{
\begin{array}{ll}
0 			& (t \le 0) 	\\
\gamma^2	& (t=\gamma^2)	\\
1-\gamma^2	& (t=\gamma - \gamma^2) \\
1			& (t \ge \gamma),
\end{array}
\right.
\]
$h$ is convex on $[0,\gamma^2]$, linear on $[\gamma^2,\gamma-\gamma^2]$, and concave on 
$[\gamma-\gamma^2,\gamma]$.  For any interval $(p,q] \subset \R$, the test function 
$r_\gamma \colon \R \to [0,1]$ defined by
\[
r_\gamma (t) ~=~
\left\{
\begin{array}{ll}
h_\gamma(t-p)		& (t \le p+\gamma^2 )	\\
1					& (p+\gamma^2 \le t \le q) \\
1 - h_\gamma(t-q)	& (t \ge q)
\end{array}
\right.
\]
converges pointwise to the characteristic function $\chi_{(p,q]}$ pointwise as $\gamma \downarrow 0$.  By  Dominated Convergence we deduce 
\[
\int r_\gamma\,d(D^2 f) ~\to~ \int \chi_{(p,q]}\, d(D^2 f) ~=~ (D^2 f)(p,q].
\]
But the left-hand side is
\begin{eqnarray*}
- \int r'_\gamma f'
& = & 
- \int_p^{p+\gamma} \Big(\frac{1}{\gamma} + O(1)\Big) f' ~-~ 
  \int_q^{q+\gamma} \Big(-\frac{1}{\gamma} + O(1)\Big) f' \\
&=& 
\frac{f(q+\gamma) - f(q)}{\gamma} ~-~ \frac{f(p+\gamma) - f(p)}{\gamma} ~+~ O(\gamma) \\
& = &
f'_+(q) - f'_+(p) + O(\gamma)
\end{eqnarray*}
as $\gamma \downarrow 0$.  We thus reproduce our earlier definition:
\[
(D^2 f)(p,q] ~=~ f'_+(q) - f'_+(p).
\]
Clearly, this fact also holds for any difference-of-convex function $f \colon \R \to \R$.
}
\end{exa}

The modulus of nonconvexity for a function $f \colon \Rn \to \R$, which we denoted $\Lambda(\delta)$ (for $\delta > 0$) is the supremum of the concave deviation of the restriction of $f$ to line segments of the form
\[
S ~=~ \{ z+tw : 0 \le t \le \delta \}
\]
for some point $z \in \Rn$ and unit direction $w \in \Rn$. 
That concave deviation is the measure of $S$ under the negative part of the measure $D^2 (f|_S)$.
We would therefore like to compare the distributional second derivative of this restriction with the distributional second derivative $D^2 f$.  As we see in the next result, we should focus specifically on the {\em directional distributional second derivative} $w^T (D^2 f)w$.

A simple approach is furnished by mollification.  We fix a {\em mollifier} $\phi \colon \Rn \to \R$:  a test function satisfying $\int\phi = 1$ and with the property that, as $\gamma \downarrow 0$, the function 
$\phi_{\gamma}(x) = \gamma^{-n}\phi(\frac{1}{\gamma} x)$ converges as a distribution to the Dirac delta function.  Given a Radon measure $\mu$ on $\Rn$ and any test function $g \colon \Rn \to \R$, we can define the convolution $g \star \mu \colon \Rn \to \R$ by
\[
(g \star \mu)(y) ~=~ \int g(y-x)\,d\mu(x).
\]

\begin{thm}[Chain rule via mollification]
Consider a locally Lipschitz function $f \colon \Rn \to \R$, a mollifier $\phi \colon \Rn \to \R$, and a direction $w \in \Rn$.  Then for almost all points $z \in \Rn$, the function $f$ is differentiable almost everywhere on the line $z+\R w$, and the distributional second derivative $D^2 h$ of the function defined by
\[
h(t) = f(z + tw) \qquad (t \in \R)
\]
is the distributional limit, as $\gamma \downarrow 0$, of the convolution
\[
t ~\mapsto~ \Big(\phi_\gamma \star \big(w^T(D^2 f)w\big) \Big)(z+tw).
\]
\end{thm}

\pf
By Rademacher's Theorem and standard properties of convolutions
\cite[Section 4.2, Theorem 1(iv)]{evans-gariepy}, there exists a full measure set $\Omega \subset \Rn$ on which $f$ is differentiable and the convolution $\phi_\gamma \star (w^T \nabla f)$ converges pointwise to the essentially bounded function $w^T\nabla f$.  By Fubini's Theorem, for almost all points $z \in \Rn$, we have $z+tw \in \Omega$ for almost all $t \in \R$.  Restricting attention to such $z$, consider any test function $g \colon \R \to \R$.  Fubini's Theorem implies
\begin{eqnarray*}
\lefteqn{
\int_{t \in \R} g(t) \int_{x \in \Rn} \phi_\gamma(z+tw-x)\, d\big(w^T(D^2 f)w\big)(x)\, dt}
\hspace{3cm} \mbox{}
\\
& & =  \int_{x} \Big( \int_{t} g(t) \phi_\gamma(z+tw-x)\, dt \Big) \, d\big(w^T(D^2 f)w\big)(x)
\\
& & = \int_{x} w^T \nabla f(x) \Big( w^T\int_{t} g(t) \nabla \phi_\gamma(z+tw-x)\, dt \Big) \, dx.
\end{eqnarray*}
(We can interchange the order of differentiation and integration since the test functions $g$ and $\phi_\gamma$ are well behaved.)  Rewriting, integrating by parts, and using Fubini's Theorem and Dominated Convergence again, we obtain
\begin{eqnarray*}
\lefteqn{
\int_{x} w^T \nabla f(x) \int_{t} g(t) w^T \nabla \phi_\gamma(z+tw-x)\, dt \, dx}
\hspace{3cm} \mbox{} 
\\
& & = 
\int_{x } w^T \nabla f(x) \int_{t} g(t) \frac{d}{dt} \phi_\gamma(z+tw-x)\, dt \, dx
\\
& & = 
 - \int_{x } w^T \nabla f(x) \int_{t} g'(t) \phi_\gamma(z+tw-x)\, dt \, dx
 \\
& & = 
 - \int_{t} g'(t) \int_{x }\phi_\gamma(z+tw-x)  w^T \nabla f(x) \, dx\, dt \\
& & = 
 - \int_{t} g'(t) \big(\phi_\gamma \star (w^T \nabla f)\big)(z+tw)\, dt \\
& & \to~
- \int_{t} g'(t) w^T \nabla f(z+tw)\, dt ~=~ - \int g' h' ~=~ \int g\,d(D^2 h),
\end{eqnarray*}
as desired.
\finpf

\noindent
In this result, the effect of the convolution is to focus attention on the line through the point $z$ in the direction $w$.  In informal language, we deduce that {\em the modulus of nonconvexity $\Lambda(\delta)$ is determined by the concentration of the negative parts of the measures $w^T (D^2 f)w$ around line segments of length $\delta$ in unit directions $w$}.

For more intuition on directional distributional second derivatives of the form $w^T(D^2f)w$, let us consider a convex function $f \colon \R^2 \to \R$.  After a suitable choice of basis we can suppose that 
$w$ is the first unit vector $e^1$ and therefore consider the Radon measure $(Df)_{11}$.  To understand this measure, consider the integral
\[
\int_{\R^2} r_\gamma(x_1) g(x_2)\, d(D^2 f)_{11}(x)
\]
for the function $r_\gamma$ of Example \ref{univariate} and any test function $g \colon \R \to \R$.  As $\gamma \downarrow 0$, we observe
\begin{eqnarray*}
\int_{\R^2} r_\gamma(x_1) g(x_2)\, d\big(e_1^T(D^2 f)e_1\big)(x)
&=&
-\int_{\R^2} \big(e_1^T \nabla f(x)\big) \cdot \Big(e_1^T \nabla\big(r_\gamma(x_1) g(x_2)\big)\Big)\,dx
\\
&=&
-\int_{\R^2} \frac{\partial f}{\partial x_1} r'_\gamma(x_1) g(x_2)\,dx_1\,dx_2 \\
& \to &
\int 
\Big(
f'(\left[\begin{smallmatrix} q \\ t \end{smallmatrix}\right] ; 
\left[\begin{smallmatrix} 1 \\ 0 \end{smallmatrix}\right]) - 
f'(\left[\begin{smallmatrix} p \\ t \end{smallmatrix}\right] ; 
\left[\begin{smallmatrix} 1 \\ 0 \end{smallmatrix}\right])
\Big)
 g(t)\,dt.
\end{eqnarray*}
More generally, this argument suggests, loosely, that $w^T (D^2 f)w$ measures the variation of the directional derivative $f'(\cdot;w)$ along the direction $w$.

\bibliographystyle{plain}
\small
\parsep 0pt

\def\cprime{$'$} \def\cprime{$'$}

\end{document}